\documentclass{amsart}
\pdfoutput=1

\usepackage{microtype}
\usepackage[mathscr]{eucal}
\usepackage{mathtools}
\usepackage{amsthm}
\usepackage{enumitem}
\usepackage{hyperref}
\usepackage{mleftright}

\DeclareMathOperator{\inte}{int}
\DeclareMathOperator{\conv}{conv}
\DeclareMathOperator{\Imm}{Imm}
\DeclareMathOperator{\vol}{vol}
\DeclareMathOperator{\diam}{diam}

\theoremstyle{plain}
\newtheorem{theorem}{Theorem}[section]

\newtheorem{lemma}[theorem]{Lemma}
\newtheorem{corollary}[theorem]{Corollary}
\theoremstyle{remark}
\newtheorem{remark}[theorem]{Remark}
\theoremstyle{definition}
\newtheorem{problem}[theorem]{Problem}

\begin{document}
	
\title[The convex hull of a convex space curve]{The convex hull of a convex space curve\\with four vertices}

\author{Jakob Bohr}
\address{Department of Applied Mathematics and Computer Science, Technical University of Denmark, Matematiktorvet, Building 303B, 2800 Kongens Lyngby, Denmark}

\author{Steen Markvorsen}
\address{Department of Applied Mathematics and Computer Science, Technical University of Denmark, Matematiktorvet, Building 303B, 2800 Kongens Lyngby, Denmark}
\email{stema@dtu.dk}

\author{Matteo Raffaelli}
\address{School of Mathematics, Georgia Institute of Technology, Atlanta, Georgia 30332}
\email{raffaelli@math.gatech.edu}
\date{December 4, 2025}

\begin{abstract}
We obtain an upper bound for the volume of the convex hull of a simple closed Frenet curve with exactly four vertices, i.e., four points of vanishing torsion, and lying on the boundary of its convex hull. Moreover, we show that the upper bound is attained when the curve intersects every plane in at most four points, a condition studied by Scherk and Segre in the 1930s. The proof relies on the fact that, under the four-vertex assumption, the convex hull is a union of line segments and therefore admits an elementary parametrization. We also comment on a question posed by Newson in 1899.
\end{abstract}
\maketitle

\section{Introduction}

In 1934, Bonnesen and Fenchel~\cite[p.~117]{bonnesen1987} proposed the following isoperimetric problem; see also \cite[section~A28]{croft1991}.

\begin{problem}\label{IsoperimetricProblem}
Among all closed curves in $\mathbb{R}^{3}$ of a given length, find the one whose convex hull has maximal volume.
\end{problem}

Problem~\ref{IsoperimetricProblem} extends naturally to curves in $\mathbb R^d$ and also has an analogue for open curves. Several special cases have been settled: the case of open convex curves in $\mathbb R^d$~\cite{egervary1949, krein1977, nudelman1975}, where ``convex'' means that no plane intersects the curve in more than three points; the case of closed curves in $\mathbb{R}^{2d}$~\cite{schoenberg1954}, again assuming that no hyperplane intersects the curve in more than $2d$ points; and the case of closed space curves subject to a symmetry constraint~\cite{melzak1960}. For additional related results, see \cite{derry1956, zalgaller1997, dedios2023}.

Although Problem~\ref{IsoperimetricProblem} has inspired many significant results about convex hulls and their volumes, most of these apply only under the restrictive convexity assumptions described above. On the other hand, closed convex curves do not exist in odd dimensions~\cite[Corollary~4.6]{amendola2023}. Consequently, little is known about the volume of the convex hull of a general closed space curve; the only relevant estimate to date appears to be the one provided by Tilli in \cite[Corollary~1.4]{tilli2010}.

The purpose of this short note is to present a volume estimate for the convex hull of a simple closed \emph{convex} curve in $\mathbb{R}^{3}$ with (positive curvature and) exactly four vertices, that is, four points where the torsion vanishes; here, and in the rest of the paper, \textit{convex} means that the curve lies on the boundary of its convex hull. Although these hypotheses are restrictive, it has been known since the work of Scherk and Segre in the 1930s that if a simple closed convex curve intersects every plane in at most four points, then it has exactly four vertices~\cite{scherk1936, segre1936}; see also \cite[p.~357]{pohl1966}. Moreover, the converse fails~\cite{ghomi2024}. Recall that four is the minimum number of vertices such a curve can have, and that the torsion changes sign at each of them~\cite{sedykh1994, ghomi2019}.

To state our main results, let $\Imm^{k\geq 1}([0,L], \mathbb{R}^{3})$ be the space of $\mathcal{C}^{k}$ curves $\beta \colon [0,L] \to \mathbb{R}^{3}$ with positive speed. We say that $\beta$ is \textit{closed} if $\beta^{(h)}(0) =\beta^{(h)}(L)$ for $h\leq k$, and that $\beta$ is \textit{simple} if it is injective on $[0,L)$. Once and for all, let $\gamma \in \Imm^{3}([0,L], \mathbb{R}^{3})$ be a simple closed convex curve with unit speed and nowhere vanishing curvature, and let $\conv(\gamma)$ be its convex hull.

\begin{theorem}\label{thm:main}
If $\gamma$ has exactly four vertices, then
\begin{equation}\label{eq:volume}
\vol(\conv(\gamma)) \leq \frac{1}{24}\int_{0}^{L}\int_{0}^{L} \lvert \det(\gamma'(t_{1}), \gamma'(t_{2}), \gamma(t_{1}) - \gamma(t_{2}))\rvert\, dt_{1}\, dt_{2},
\end{equation}
where $\vol$ denotes volume.
\end{theorem}

Clearly, the integrand in \eqref{eq:volume} is bounded from above by $\diam(\gamma) \coloneqq \max_{[0,L]^{2}}\lVert \gamma(t_{1})-\gamma(t_{2})\rVert$. Hence the above theorem immediately yields the following inequality.

\begin{corollary}
If $\gamma$ has exactly four vertices, then
\begin{equation*}
\vol(\conv(\gamma)) < \frac{1}{24} \diam(\gamma)L^{2}.
\end{equation*}
\end{corollary}

\begin{remark}
Since $\gamma$ is closed, we have $L > 2\diam(\gamma)$, which implies $\vol(\conv(\gamma)) < L^{3}/48$ when the number of vertices is four. This may be compared with the general estimate in \cite[Corollary~1.3]{tilli2010}, which gives $\vol(\conv(\gamma)) < L^{3}/36$.
\end{remark}

\begin{remark}
If $\gamma$ lies on a sphere of radius $R$, then $\diam(\gamma) \leq 2R$. So $\vol(\conv(\gamma)) < R L^{2}/12$ when it has exactly four vertices.
\end{remark}

The proof of Theorem~\ref{thm:main}, whose details are presented in the next two sections, may be summarized as follows. When $\gamma$ has exactly four vertices, its convex hull is a union of \emph{chords}, i.e., line segments with endpoints on $\gamma([0,L])$, and thus admits a parametrization of the form $(t_{1},t_{2},u)\mapsto \gamma(t_{1})+u(\gamma(t_{2}) - \gamma(t_{1}))$. Besides, the convexity assumption ensures that $\gamma$ is star-shaped with respect to any point $p$ in the interior of its convex hull. By projecting $\gamma$ radially onto a sphere centered at $p$, we obtain a spherical curve $\gamma_{p}$ that has as many pairs of antipodal points as there are chords passing through $p$. In particular, the projected curve must have at least two pairs of antipodal points; otherwise, $p$ would lie on the boundary of the convex hull. This implies that $p$ is covered at least four times by the parametrization, and the desired inequality follows.

In \cite{segre1936}, Segre proves that if $\gamma$ meets every plane in at most four points (counted with multiplicity), then not only is the convex hull a union of chords, but any point in its interior is contained in exactly two chords. Hence our proof of Theorem~\ref{thm:main} leads to the following corollary.

\begin{corollary}\label{cor:formula}
If no plane intersects $\gamma$ in more than four points (counted with multiplicity), then 
\begin{equation}\label{eq:formula}
\vol(\conv(\gamma)) = \frac{1}{24}\int_{0}^{L}\int_{0}^{L} \lvert \det(\gamma'(t_{1}), \gamma'(t_{2}), \gamma(t_{1}) - \gamma(t_{2}))\rvert\, dt_{1}\, dt_{2}.
\end{equation}
\end{corollary}

As we explain in section~\ref{sec:newson}, Corollary~\ref{cor:formula} gives a partial answer to an old question of Newson~\cite{newson1899}.

\section{The radial projection}

Let $o$ be the origin of $\mathbb{R}^{3}$, and suppose that $o \in \inte \conv(\gamma)$. Note that since $\gamma$ is convex, no ray emanating from $p \in \inte\conv(\gamma)$ intersects $\gamma$ in more than one point. Hence the \textit{radial projection (of $\gamma$ with respect to $p$)} 
\begin{equation*}
\gamma_{p} \coloneqq (\gamma-p)/\lVert \gamma-p\rVert
\end{equation*}
is a well-defined simple closed curve on $\mathbb{S}^{2}$, and it divides $\mathbb S^2$ into two connected regions. The purpose of this section is to record the following basic observation, which will play a key role in the proof of Theorem~\ref{thm:main}.

\begin{lemma}\label{lm:antipodal}
Suppose that $\gamma_{o}\coloneqq \gamma/\lVert\gamma\rVert$ intersects its antipodal reflection at a single pair of points $\pm q \in \mathbb{S}^{2}$, and let $n \in T_q\mathbb S^2$ be the unit normal of $\gamma$ pointing to the spherical region of smaller area. Then, for any sufficiently small $p\in \inte \conv(\gamma)$ satisfying $\langle p, n\rangle < 0$, the curve $\gamma_{p}$ is disjoint from its antipodal reflection.
\end{lemma}

\begin{proof}
To begin with, note that when $\gamma_{o}$ and its antipodal reflection intersect at a single pair of points, they must do so tangentially, as otherwise they would cross each other~\cite[Note~8.2]{ghomi2013}; hence $n$ is also normal to $\gamma_o$ at $-q$. 

Let $u\coloneqq p/\lVert p\rVert$, and let $F_t\colon r\mapsto (\gamma(t)-r u)/\lVert \gamma(t)-r u\rVert$. A computation shows that
\begin{equation*}
F_t'(0) = -\frac{u}{\lVert \gamma(t)\rVert} + \frac{\langle\gamma(t),u\rangle}{\lVert \gamma(t)\rVert^3}\gamma(t).
\end{equation*} 
In particular, if $t_+,t_-$ satisfy $\gamma_o(t_+)=q$ and $\gamma_o(t_-)=-q$, then
\begin{equation*}
\langle F_{t_\pm}'(0),n \rangle= -\frac{\langle u,n \rangle}{\lVert\gamma(t_\pm)\rVert}.
\end{equation*} 
The above equation shows that the infinitesimal motion produced by $p$ points to the same region of $\mathbb S^2$ at $q$ and $-q$. 

Next, pick geodesic normal coordinates at $q$ whose axes are the unit tangent $\gamma_o'(t_+)$ and the normal $n$. In these coordinates, the curves $\pm\gamma_p$ may be represented by functions $y = f_\pm(r,x)$ around $t_\pm$, as follows. Write $\gamma_p(t)=(X_+(r,t),Y_+(r,t))$, and let $T_+(r,\cdot)$ be the inverse of $X_+(r,\cdot)$, which exists by the inverse function theorem. Then
\begin{equation*}
f_+(r,x) = Y_+(r, T_+(r,x)).
\end{equation*} 
Similarly,
\begin{equation*}
f_-(r,x) = Y_-(r, T_-(r,x)).
\end{equation*} 
Our problem is equivalent to showing that for all sufficiently small $r>0$, the function $x\mapsto h(r,x)\coloneqq f_+(r,x)-f_-(r,x)$ does not intersect the $x$-axis in a neighborhood of $0$.

Note that $h(0,x) \geq 0$, because $\gamma_o$ and $-\gamma_o$ cannot cross. Moreover, since $\partial_x f_\pm(0,0)=0$, we have
\begin{equation*}
\partial_rf_\pm(0,0)=\partial_r Y_\pm(0,t_\pm)=\pm \langle F_{t_\pm}'(0),n \rangle,
\end{equation*} 
which implies
\begin{equation*}
\partial_r h(0,0)= -\frac{\langle u,n \rangle}{\lVert\gamma(t_+)\rVert}-\frac{\langle u,n \rangle}{\lVert\gamma(t_-)\rVert} >0.
\end{equation*} 
By continuity of $\partial_r h$, there exist $\varepsilon >0$ and $\delta >0$ such that $\partial_r h(r,x)>0$ when $\lvert r\rvert \leq \varepsilon$ and $\lvert x\rvert\leq\delta$. Now fix any $0<r\leq \varepsilon$. The fundamental theorem of calculus gives
\begin{equation*}
h(r,x) = h(0,x)+\int_0^r \partial_r h(s,x)\, ds >0.
\end{equation*} 
Hence, for all sufficiently small $x$ and $r>0$, no intersections between $f_+$ and $f_-$ occur, and the lemma follows. 
\end{proof}

\section{Proof of Theorem~\ref{thm:main}}

Here we prove Theorem~\ref{thm:main}. The first step is to show that under the four-vertex assumption, any point in the convex hull of $\gamma$ can be written as a convex combination of just two points of $\gamma([0,L])$.

\begin{lemma}[\cite{ivanishvili2015, zatiskiy2016}]\label{lm:chords}
If $\gamma$ has exactly four vertices, then $\conv(\gamma)$ is a union of \textit{chords}, i.e., line segments with endpoints in $\gamma([0,L])$.
\end{lemma}

\begin{proof}
Since $\gamma$ is connected, a theorem of Fenchel~\cite{fenchel1929, hanner1951} tells us that every point in $\conv(\gamma)$ can be expressed as a convex combination of at most three points of $\gamma$. Under the hypothesis of the lemma, \cite[Corollary~2]{romerofuster1997} then shows that the boundary of $\conv(\gamma)$ is a union of chords. Consequently, as explained in \cite{ivanishvili2015}, the interior has the same property. The argument is reproduced below for the reader’s convenience.
 
Let $p$ be a point in the interior, which we may assume to coincide with the origin of $\mathbb{R}^{3}$, and let $p_{1}$ and $p_{2}$ be the intersections of the $z$-axis with the boundary of $\conv(\gamma)$. Then we know by assumption that $p_{i}$ ($i=1,2$) belongs to a line segment with endpoints $q_{i},r_{i}\in \gamma([0,L])$. Next, extending $q_{1},q_{2}$ to a continuous map $q\colon[1,2] \to\gamma([0,L])$, define $r_{t}$ as the intersection of the image of $\gamma$ with the plane spanned by $q_{t}$ and the $z$-axis, and $p_{t}$ as the intersection of the chord $(q_{t},r_{t})$ with the $z$-axis. It is clear that $p_{t}$ depends continuously on $t$ and agrees with $p_{1}, p_{2}$ for $t=1,2$. Hence there exists $t\in[1,2]$ such that the chord $(q_{t},r_{t})$ intersects the $z$-axis in the point $p$, which is the desired conclusion.
\end{proof}

Next we show that any point in the interior of the convex hull is hit by at least two chords.

\begin{lemma}\label{lm:boundary}
Suppose that $\conv(\gamma)$ is a union of chords. If $p \in \conv(\gamma)$ belongs to only one chord, then $p \in \partial\conv(\gamma)$.
\end{lemma}

\begin{proof}
By contradiction, suppose that $p_0\in\inte\conv(\gamma)$ lies on exactly one chord. Observe that for any $p\in\inte\conv(\gamma)$, the radial projection $\gamma_p$ has as many pairs of antipodal points as there are chords through $p$. After translating so that $p_0=o$, Lemma~\ref{lm:antipodal} implies that for every sufficiently small $p$ with $\langle p, n\rangle <0$, the curves $\gamma_p$ and $-\gamma_p$ are disjoint. Since $\conv(\gamma)$ is a union of chords, it follows that every neighborhood of $o$ contains points outside $\conv(\gamma)$. Hence $o\in\partial\conv(\gamma)$, contradicting the assumption that $o\in\inte\conv(\gamma)$.
\end{proof}

Now we are ready to finalize our proof.

\begin{proof}[Proof of Theorem~\textup{\ref{thm:main}}]
By Lemma~\ref{lm:chords}, the map $\sigma \colon [0,L]^{2}\times [0,1] \to \mathbb{R}^{3}$ defined by
\begin{equation*}
\sigma(t_{1}, t_{2},u)=\gamma(t_{1})+u(\gamma(t_{2})-\gamma(t_{1}))
\end{equation*}
is a parametrization of $\conv(\gamma)$. In particular, note that $\conv(\gamma)$ is covered at least twice by $\sigma$, because $\sigma(t_{1},t_{2},u)=\sigma(t_{2},t_{1},1-u)$. Furthermore, Lemma~\ref{lm:boundary} implies that any point in the interior of $\conv(\gamma)$ is contained in at least two chords of $\gamma$. Hence the volume of $\conv(\gamma)$ satisfies
\begin{equation*}
\vol(\conv(\gamma))\leq \frac{1}{4}\iiint_{[0,L]^{2}\times [0,1]}\lvert \det(J_{\sigma}) \rvert\, dV,
\end{equation*}
where $J_{\sigma}$ denotes the Jacobian matrix of $\sigma$. Since
\begin{equation*}
\det(J_{\sigma}) = \mleft(u-u^{2}\mright)\det(\gamma'(t_{1}), \gamma'(t_{2}), \gamma(t_{2})-\gamma(t_{1})),
\end{equation*}
integration with respect to $u$ gives
\begin{equation*}
\vol(\conv(\gamma))\leq \frac{1}{24} \int_{0}^{L}\int_{0}^{L}\lvert \det(\gamma'(t_{1}), \gamma'(t_{2}), \gamma(t_{2})-\gamma(t_{1})) \rvert\, dt_{1}\,dt_{2},
\end{equation*}
as desired.
\end{proof}

\section{Newson's challenge}\label{sec:newson}
Here we briefly comment on a question posed by Newson in 1899.

Let $\beta \in \Imm^{1}([0,L], \mathbb{R}^{3})$ be a simple closed unit-speed curve lying on a plane containing the origin. Then the area enclosed by $\beta$ is given by
\begin{equation}\label{eq:area}
\frac{1}{2} \biggl\lVert \int_{0}^{L}\beta'(t)\times \beta(t) \, ds \biggr\rVert.
\end{equation}
In \cite{newson1899}, Newson derived the area formula (or, rather, an equivalent version of it) by interpreting $\beta([0,L])$ as a polygon with infinitely many sides, i.e., as the limiting case of the expression
\begin{equation*}
\frac{1}{2} \biggl\lVert \sum_{i=1}^{n} \beta_{i} \times \beta_{i+1}\biggr\rVert,
\end{equation*}
which measures the area of an inscribed polygon with $n$ vertices. He then challenged the reader of the Annals to come up with a three-dimensional version of \eqref{eq:area} by interpreting a closed surface as a polyhedron with infinitely many faces. While we do not provide a complete answer to this problem, in the following we explain how Newson's idea may be used to derive formula~\eqref{eq:formula}.

As before, let $\gamma \in \Imm^{3}([0,L], \mathbb{R}^{3})$ be a simple closed unit-speed curve with nowhere vanishing curvature, lying on the boundary of its convex hull, and suppose that no plane intersects $\gamma$ in more than four points. Partitioning $[0,L]$ into $n-1$ sufficiently small subintervals $[t_{i}, t_{i+1}]$, we may approximate $\gamma$ by the polygonal line with vertices $\gamma_{i}\coloneqq \gamma(t_{i})$. Then, as $n$ grows large,
\begin{equation*}
\bigcup_{i,j} T_{ij} \to \conv(\gamma),
\end{equation*}
where $T_{ij}$ denotes the tetrahedron $(\gamma_{i}, \gamma_{i+1}, \gamma_{j},\gamma_{j+1})$, and each point in the interior of $T_{ij}$ lies in the interior of exactly two tetrahedra. Since $T_{ij}$ has signed volume
\begin{equation*}
\frac{1}{6}\det (\gamma_{i+1}-\gamma_{i}, \gamma_{j}-\gamma_{i}, \gamma_{j+1}-\gamma_{j}),
\end{equation*}
it follows that 
\begin{equation*}
\vol(\conv(\gamma_i)) \approx \frac{1}{24}\sum_{i,j} \lvert\det (\gamma_{i+1}-\gamma_{i}, \gamma_{j}-\gamma_{i}, \gamma_{j+1}-\gamma_{j}) \rvert,
\end{equation*}
where the extra factor of $1/2$ comes from the equality $T_{ij}=T_{ji}$. Thus, letting $n\to \infty$, we obtain \eqref{eq:formula}.

\section*{Acknowledgments}
We thank Carlos Am\'{e}ndola, Mohammad Ghomi, and Paata Ivanishvili for useful communications.

\bibliographystyle{amsplain}
\bibliography{biblio}
\end{document}